\documentclass[12pt]{article}

\title{Leibniz, Information, Math and Physics}

\author{G. J. Chaitin\thanks{IBM Research Division, P. O. Box 218, Yorktown Heights, NY 10598, USA,
E-mail: \emph{chaitin@us.ibm.com}}}

\date{}

\begin{document}

\maketitle

\begin{abstract}
The information-theoretic point of view proposed by Leibniz in 1686
and developed by algorithmic information theory (AIT) suggests that mathematics
and physics are not that different.
This will be a first-person account of some doubts and speculations about
the nature of mathematics that I have entertained for the past three decades,
and which have now been incorporated in a digital philosophy paradigm shift  
that is sweeping across the sciences. 
\end{abstract}

\section*{1. What is algorithmic information theory?}

The starting point for my own work on AIT forty years ago was the insight
that a scientific theory is a computer program that calculates the observations,
and that the smaller the program is, the better the theory.  If there is no theory,
that is to say, no program substantially smaller than the data itself,
considering them both to be finite binary strings, then
the observations are algorithmically random, theory-less, unstructured, incomprehensible
and irreducible.

\begin{center}
   theory = program $\longrightarrow$ \fbox{\textbf{Computer}} 
   $\longrightarrow$ output = experimental data  
\end{center}

So this led me to a theory of randomness based on program-size complexity [1], whose
main application turned out to be not in science, but in mathematics, more specifically,
in meta-mathematics, where it yields powerful new information-theoretic versions of G\"odel's
incompleteness theorem [2, 3, 4]. (I'll discuss this in Section 3.)

And from this new information-theoretic point of view, math and physics do not
seem too different. In both cases understanding is compression, and is measured by
the extent to which empirical data and mathematical theorems are respectively compressed into
concise physical laws or mathematical axioms, both of which are embodied in computer software
[5].

And why should one use reasoning at all in mathematics?!  Why not proceed entirely
empirically, more or less as physicists do?  Well, the advantange of proving things is that 
assuming a few bits of axioms is less risky than assuming many empirically-suggested 
mathematical assertions. (The disadvantage, of course, is the length of the proofs
and the risk of faulty proofs.) Each bit in an irreducible axiom of a mathematical theory is
a freely-chosen independent assumption, with an \emph{a priori} probability 
of half of being the right choice,
so one wants to reduce the number of such independent choices to a minimum in
creating a new theory.  

So this point of view would seem to suggest that while math and physics are admittedly
different, perhaps they are not as different as most people usually believe. Perhaps we should
feel free to pursue not only rigorous, formal modern proofs, but also the swash-buckling
experimental math that Euler enjoyed so much.  And in fact theoretical computer
scientists have to some extent already done this, 
since their $\mathbf{P} \neq \mathbf{NP}$ hypothesis is probably currently the best
candidate for canonization as a new axiom. 
And, as is suggested in [6], another possible candidate is
the Riemann hypothesis.

But before discussing this in more detail, 
I'd like to tell how I discovered that in 1686 Leibniz anticipated some of the
basic ideas of AIT.

\section*{2. How Leibniz almost invented algorithmic information theory [7]}

One day last year, while preparing my first philosophy paper [5], for a philosophy congress in Bonn,
I was reading a little book on philosophy by Hermann Weyl
that was published in 1932, and I was amazed to
find the following, which captures the essential idea of my definition of algorithmic randomness:

\begin{quote}
``The assertion that nature is governed by strict laws
is devoid of all content if we do not add the statement that it is governed
by mathematically simple laws\ldots\ That {\bf the notion of law becomes empty
when an arbitrary complication is permitted} was already pointed out
by Leibniz in his {\it Metaphysical Treatise\/} [{\it Discourse on Metaphysics\/}].
Thus simplicity becomes a working principle in the natural sciences.

---Weyl [8, pp.\ 40--42]. See a similar
discussion on pp.\ 190--191 of Weyl [9], Section 23A, ``Causality and Law''.
\end{quote}

In fact, I actually read Weyl [9] as a teenager, 
before inventing AIT at age 15,
but the matter is not stated so
sharply there.  And a few years ago I stumbled on the above-quoted text in Weyl [8],
but hadn't had the time to pursue it until stimulated to do so
by an invitation from the German Philosophy Association to talk at their 2002 annual congress,
that happened to be on limits and how to transcend them.

So I got a hold of Leibniz's \emph{Discourse on Metaphysics} 
to see what he actually said.  Here it is:

\begin{quote}
``As for the simplicity of the ways of God, this holds
properly with respect to his means, as opposed to the variety, richness,
and abundance, which holds with respect to his ends or effects.''

``\ldots not only does nothing completely irregular occur in the world, but
we would not even be able to imagine such a thing. Thus, let us assume, for
example, that someone jots down a number of points at random on a piece of
paper, as do those who practice the ridiculous art of geomancy.\footnote
{[A way to foretell the future; a form of divination.]}
I maintain that it is possible to find a geometric line whose [m]otion is constant
and uniform, following a certain rule, such that this line passes through all the
points in the same order in which the hand jotted them down.''

``But, {\bf when a rule is extremely complex, what is in
conformity with it passes for irregular}. Thus, one can say, in whatever
manner God might have created the world, it would always have been regular
and in accordance with a certain general order. But {\bf God has chosen the
most perfect world, that is, the one which is at the same time the simplest
in hypotheses and the richest in phenomena}, as might be a line in geometry
whose construction is easy and whose properties and effects are extremely
remarkable and widespread.''

---Leibniz, {\it Discourse on Metaphysics,} 1686, Sections 5--6, 
as translated by Ariew and Garber [10, pp.\ 38--39].
\end{quote}

\begin{center}
   ideas = input $\longrightarrow$ \fbox{\textbf{Mind of God}} 
   $\longrightarrow$ output = the universe
\end{center}

And after finishing my paper [5] for the Bonn philosophy congress,
I learned that Leibniz's original \emph{Discourse on Metaphysics} was in French, 
which I know, and fortunately not in Latin, which I don't know,
and that it was readily available from France:

\begin{quote}
``Pour ce qui est de la simplicit\'e des voyes de Dieu, elle a lieu proprement \`a
l'\'egard des moyens, comme au contraire la variet\'e, richesse ou abondance y a lieu
\`a l'\'egard des fins ou effects.''

``\ldots non seulement rien n'arrive dans le monde, qui soit absolument irregulier, mais
on ne s\c{c}auroit m\^emes rien feindre de tel. Car supposons par exemple que quelcun
fasse quantit\'e de points sur le papier \`a tout hazard, comme font ceux qui exercent
l'art ridicule de la Geomance, je dis qu'il est possible de trouver une ligne
geometrique dont la [m]otion soit constante et uniforme suivant une certaine regle, en
sorte que cette ligne passe par tous ces points, et dans le m\^eme ordre que la main
les avoit marqu\'es.''

``Mais \textbf{quand une regle est
fort compos\'ee, ce qui luy est conforme, passe pour irr\'egulier}. Ainsi on peut
dire que de quelque maniere que Dieu auroit cr\'e\'e le monde, il auroit
tousjours est\'e regulier et dans un certain ordre general. Mais \textbf{Dieu a choisi
celuy qui est} le plus parfait, c'est \`a dire celuy qui est en m\^eme temps \textbf{le plus
simple en hypotheses et le plus riche en phenomenes,}
comme pourroit estre une ligne de Geometrie dont la construction
seroit ais\'ee et les propriet\'es et effects seroient
fort admirables et d'une grande \'etendue.''

---Leibniz, \emph{Discours de m\'etaphysique,} \textbf{V--VI} [11, pp.\ 40--41].
\end{quote}

(Here ``dont la motion'' is my correction. The Gallimard text [11] states ``dont la notion,'' an
obvious misprint, which I've also corrected in the English translation by Ariew and Garber.)

So, in summary, Leibniz observes that for any finite set of points there
is a mathematical formula that produces a curve that goes through them all, and it can be parametrized
so that it passes through the points in the order that they were given and with a constant
speed.  So this cannot give us a definition of what it means for a set of points to obey
a law.  But if the formula is very simple, and the data is very complex, 
then that's a \emph{real law}!

Recall that Leibniz was at the beginning of the modern era, in which ancient metaphysics
was colliding with modern empirical science.  And he was a great mathematician as well
as a philosopher.  So here he is able to take a stab at clarifying what it means to say
that Nature is lawful and what are the conditions for empirical science to be possible.

AIT puts more meat on Leibniz's proposal, it makes his ideas more precise
by giving a precise definition of complexity.

And
AIT goes beyond Leibniz by using program-size complexity to clarify what it means for a
sequence of observations to be lawless, one which has no theory, and by applying this to studying the
limits of formal axiomatic reasoning, i.e., what can be achieved by mindlessly and mechanically
grinding away deducing all possible consequences of a fixed set of axioms.
(I'll say more about metamathematical applications of AIT in Section 3 below.)

\begin{center}
   axioms = program $\longrightarrow$ \fbox{\textbf{Computer}} 
   $\longrightarrow$ output = theorems
\end{center}

By the way, the articles by philosophy professors that I've seen that discuss the above text by Leibniz
criticize what they see as the confused and ambiguous nature of his remarks.  
On the contrary, I admire
his prescience and the manner in which he has unerringly identified the central issue, the key
idea. He even built a mechanical calculator and with his 
speculations regarding a \emph{Characteristica Universalis} 
(``Adamic'' language of creation)
envisioned something that Martin Davis [12] has argued was a direct intellectual ancestor of
the universal Turing machine, which is precisely the device that 
is needed in order for AIT to be able to quantify Leibniz's original insight!

Davis quotes some interesting remarks by Leibniz about the practical utility of his calculating
machine. Here is part of the Davis Leibniz quote:
\begin{quote}
``And now that we may give final praise to the machine we may say that it will be desirable
to all who are engaged in computations which, it is well known, are the managers of
financial affairs, the administrators of others' estates, merchants, surveyors, geographers,
navigators, astronomers\ldots\  For it is unworthy of excellent men to lose hours like
slaves in the labor of calculations which could safely be relegated to anyone else if the
machine were used.''
\end{quote}

This reminds me of a transcript of a lecture that von Neumann 
gave at the inauguration of the NORC (Naval Ordnance Computer)
that I read many years ago. It attempted
to convince people that computers were of value.
It was a hard sell!
The obvious 
practical and scientific
utility of calculators and computers, though it was evident to Leibniz, Babbage and von Neumann,
was far from evident to most people.
Even von Neumann's colleagues at the Princeton Institute of Advanced Study completely
failed to understand this (see Casti [13]). 

And I am almost forgetting something important that 
I read in E. T. Bell [14] as a child, which is that
\textbf{Leibniz invented base-two binary notation} for integers. 
Bell reports that
this was a result of Leibniz's interest in Chinese culture; no doubt
he got it from the \emph{I Ching}.
So in a sense, all of information theory derives from Leibniz, for he was the first to emphasize
the creative combinatorial potential of the 0 and 1 bit, and how everything can be built
up from this one elemental choice, from these two elemental possibilities.
So, perhaps not entirely seriously, I should propose changing the name of the unit of information
from the \emph{bit} to the \emph{leibniz}!

\section*{3. The halting probability $\Omega$ and information-theoretic incompleteness}

Enough philosophy, let's do some mathematics!  The first step is to pick a universal binary computer
$U$ with the property that for any other binary computer $C$ there is a binary prefix $\pi_C$ 
such that
\[
   U(\pi_C \, p) = C(p).
\]
Here $p$ is a binary program for $C$ and
the prefix $\pi_C$ tells $U$ how to \emph{simulate} $C$ and does not depend on $p$. 
In the $U$ that I've picked,
$\pi_C$ consists of a description of $C$ written in the high-level non-numerical functional
programming language LISP, which is much like a computerized version of set theory, except
that all sets are finite.

Next we define the \emph{algorithmic information content} 
(program-size complexity)
of a LISP symbolic expression (S-expression) $X$ to be the
size in bits $|p|$ of the smallest binary program $p$ that makes our chosen $U$ compute $X$:
\[
   H(X) \equiv \min_{U(p)=X} |p|.
\]

Similarly, the information content or complexity of a formal axiomatic theory with the
infinite set of theorems $T$ is defined to be the size in bits of the smallest program
that makes $U$ generate the infinite set of theorems $T$, which is a set of
S-expressions.
\[
   H(T) \equiv \min_{U(p)=T} |p|.
\]
Think of this as the minimum number of bits required to tell $U$ how to run through all
possible proofs and systematically generate all the consequences of the fixed set of axioms.
$H(T)$ is the size in bits of the most concise axioms for $T$.

Next we define the celebrated halting probability $\Omega$:
\[
   \Omega \equiv \sum_{U(p)\,\mathrm{halts}} 2^{-|p|}.
\]
A small technical detail: To get this sum to converge it is necessary that programs for $U$
be ``self-delimiting.'' I.e., no extension of a valid program is a valid program, the set of
valid programs has to be a prefix-free set of bit strings.  

So $\Omega$ is now a specific,
well-defined real number between zero and one, and let's consider its binary expansion, i.e.,
its base-two representation.  Discarding the initial decimal (or binary) point, 
that's an infinite binary sequence $b_1b_2b_3\ldots\;$ To eliminate any ambiguity in case
$\Omega$ should happen to be a dyadic rational (which it actually isn't), 
let's agree to change 1000\ldots\ to 0111\ldots\ here if necessary.

Right away we get into trouble.  From the fact that knowing the first $N$ bits of $\Omega$
\[
   \Omega_N \equiv b_1b_2b_3\ldots b_N
\]
would enable us to answer the halting problem for every program $p$ for $U$ with $|p|\leq N$,
it is easy to see that the bits of $\Omega$ are \emph{computationally irreducible}:
\[
   H(\Omega_N) \geq N-c.
\]
And from this it follows using a straight-forward program-size argument (see [3]) that
the bits of $\Omega$ are also \emph{logically irreducible.}  

What does this mean?  Well, consider
a formal axiomatic theory with theorems $T$, an infinite set of S-expressions.
If we assume that a theorem of the form 
``The $k$th bit of $\Omega$ is 0/1'' is in $T$ only if it's true,
then $T$ cannot enable us to determine more than $H(T)+c'$ bits of $\Omega$.

So the bits of $\Omega$ are irreducible mathematical facts, they are mathematical facts
that contradict Leibniz's principle of sufficient reason by \textbf{being true for no reason}.
They must, to use Kantian terminology, be apprehended as things in themselves.
They cannot be deduced as consequences of any axioms or principles that are simpler than they are.

(By the way, this also implies that the bits of $\Omega$ are statistically random, e.g.,
$\Omega$ is absolutely Borel normal in every base. I.e., all blocks of digits of the same
size have equal limiting relative frequency, regardless of the radix chosen for representing 
$\Omega$.)
 
Furthermore, in my 1987 Cambridge University Press monograph [15] I celebrate the fact
that the bits of $\Omega$ can be encoded via a diophantine equation.  There I exhibit an exponential
diophantine equation $L(k,\mathbf{x})=R(k,\mathbf{x})$ with parameter $k$ and about twenty-thousand
unknowns $\mathbf{x}$ that has infinitely many solutions iff the $k$th bit of $\Omega$ is a 1.
And recently Ord and Kieu [16] have shown that this can also be accomplished 
using the even/odd parity 
of the number of solutions, rather than its finite/infinite cardinality.
So $\Omega$'s irreducibility also infects elementary number theory!

These rather brutal incompleteness results show how badly mistaken Hilbert was to assume 
that a fixed formal axiomatic theory could encompass all of mathematics.  
And if you have to extend the foundations of mathematics by constantly adding new axioms, new 
concepts and fundamental principles, then mathematics becomes much more tentative and begins
to look much more like an empirical science. At least I think so, and you can even
find quotes by G\"odel that I think point in the same direction.

These ideas are of course controversial; see for example a highly critical review of two
of my books in the \emph{AMS Notices} [17].  I discuss the hostile reaction of the logic community
to my ideas in more detail in an interview with performance artist Marina Abramovic [18].
Here, however, I prefer to tell why I think that the world is actually moving rather quickly 
in my direction.
In fact, I believe that my ideas are now part of an unstoppable tidal wave of change  
spreading across the sciences!

\section*{4. The digital philosophy paradigm shift}

As I have argued in the second half of my 2002 paper in the \emph{EATCS Bulletin} [19],
what we are witnessing now is a dramatic convergence of mathematics with theoretical
computer science and with theoretical physics.
The participants in this paradigm shift believe 
that \emph{information} and \emph{computation} are fundamental
concepts in all three of these domains, and that what physical systems actually do is
computation, i.e., information processing.  In other words, as is asked on the cover
of a recent issue of \emph{La Recherche} with an article [20] about this, ``Is God a Computer?''

But that is not quite right. Rather, we should ask, ``Is God a Programmer?''
The intellectual legacy of the West, and 
in this connection
let me recall Pythagoras, Plato, Galileo and James Jeans,
states that ``Everything is number; God is a mathematician.''
We are now beginning to believe something slightly different, a refinement of the original
Pythagorean credo: ``Everything is software; God is a computer programmer.''
Or perhaps I should say: ``All is algorithm!'' 
Just as DNA programs living beings, God programs the universe.

In the digital philosophy  movement I would definitely include: the extremely active field
of quantum information and quantum computation [21], 
Wolfram's work [22] on \emph{A New Kind of Science}, 
Fredkin's work on reversible cellular automata and his 
website at {\small \texttt{http://digitalphilosophy.org}} 
(the pregnant phrase ``digital philosophy'' is due to
Fredkin), 
the Bekenstein-t'Hooft ``holographic principle'' [23],
and AIT. 
Ideas from theoretical physics and theoretical computer 
science are definitely leaking across the traditional boundaries between these two fields.
And this holds for AIT too, because its two central concepts are versions of 
\emph{randomness} and of \emph{entropy}, 
which are ideas that I took with me from physics and into mathematical logic.
 
Wolfram's work is particularly relevant to our discussion of the nature of mathematics,
because he believes that most simple systems are either trivial or equivalent to a universal
computer, and therefore that mathematical questions are either trivial or can never be 
solved, except, so to speak, for a set of measure zero.  This he calls his \emph{principle of
computational equivalence,} and it leads him to take the incompleteness phenomenon much more seriously
than most mathematicians do.  In line with his thesis, his book 
presents a great deal of computational evidence, but not many proofs.

Another important issue studied in Wolfram's book [22] is the question of whether, to use Leibnizian
terminology, mathematics is necessary or is contingent.  
I.e., would intelligent creatures on another planet necessarily discover the same concepts that we have, or might
they develop a perfectly viable mathematics that we would have a great deal of trouble in recognizing 
as such?  Wolfram gives a number of examples that suggest that the latter is in fact the case.

I should also mention some recent books on the quasi-empirical view of mathematics [24] and on
experimental mathematics [25, 26], as well as Douglas Robertson's two volumes [27, 28] on information as a key
historical and cultural parameter and motor of social change,
and John Maynard Smith's related books on biology [29, 30]. 

Maynard Smith and Szathm\'ary [29, 30]
measure biological evolutionary progress
in terms of abrupt improvements in the way information is represented and transmitted inside
living organisms. 
Robertson sees social evolution as driven by the same motor.
According to Robertson [27, 28],
spoken language defines the human, writing creates civilization, the printing press provoked the
Renaissance, and the Internet is weaving a new World-Wide Web.
These are abrupt improvements in the way human society is able to store and transmit information.
And they result in
abrupt increases in cultural complexity, 
in abrupt increases in social \emph{intelligence,} as it were.

(And for the latest results on $\Omega$, see Calude [31].)

\section*{5. Digital philosophy is Leibnizian; Leibniz's legacy}

None of us who made this paradigm shift happen were students of Leibniz, but he 
anticipated us all.  As I hinted in a letter to \emph{La Recherche,} in a sense all of
Wolfram's thousand-page book is the development of \emph{one sentence} in Leibniz:
\begin{quote}
``\textbf{Dieu a choisi celuy qui est\ldots\ le plus simple en hypotheses et le plus riche en
                                 phenomenes}" 
\\
{\small
[God has chosen that which is the most simple in hypotheses and the most rich in
                                  phenomena] 
}
\end{quote}
This presages Wolfram's basic insight that simple programs can have very complicated-looking output.

And all of \emph{my} work may be regarded as the development of another sentence in Leibniz:
\begin{quote}
``\textbf{Mais quand une regle est fort compos\'ee, ce qui luy est conforme, passe pour
                                  irr\'egulier}''
\\
{\small
[But when a rule is extremely complex, that which conforms to it passes for random] 
}
\end{quote}
Here I see the germ of my definition of algorithmic randomness and irreducibility.

Newtonian physics is now receding into the dark, distant intellectual past.
It's not just that it has been superseded by quantum physics.
No, it's much deeper than that.  In our new interest in \textbf{complex} systems,
the concepts of energy and matter take second place to the concepts of
information and computation.
And the continuum mathematics of Newtonian physics now takes second place to
the combinatorial mathematics of complex systems.  

As E. T. Bell stated so forcefully [32],
Newton made one big contribution to math, involving the continuum, but Leibniz made \emph{two}:
his work on the continuum and his work on discrete combinatorics (which Leibniz named).
Newton obliterated Leibniz and stole from him both his royal patron and the credit for the calculus.
Newton was buried with full honors at Westminster Abbey, while a forgotten Leibniz was accompanied
to his grave by only his secretary.  But, as E. T. Bell stated a half a century ago [32],
with every passing year, the shadow cast by Leibniz gets larger and larger.

How right Bell was!
The digital philosophy paradigm is a direct intellectual descendent of Leibniz, it is part of the
Leibnizian legacy. The human race has finally caught up with this part of Leibniz's thinking.
Are there,
Wolfram and I wonder,
more treasures there that we have not yet been able to decipher and appreciate?

\section*{6. Acknowledgment; Coda on the continuum and the Kabbalah}

The author wishes to thank Fran\c{c}oise Chaitin-Chatelin for sharing 
with him her understanding and appreciation of Leibniz,  
during innumerable lengthy conversations.   
In her opinion, however, this essay does Leibniz an injustice by completely
ignoring his deep interest in the ``labyrinth of the continuum,''
which is her specialty.

Let me address her concern.  According to Leibniz, the integers are human, the discrete
is at the level of Man.  But the continuum transcends Man and brings us closer to God.
Indeed, $\Omega$ is transcendent, and may be regarded as the concentrated essence of
mathematical creativity. In a note on the Kabbalah, which regards Man as perfectable
and evolving towards God, Leibniz [33, pp.\ 112--115]
observes that with time we shall know all interesting theorems with proofs of up to any
given fixed size, and this can be used
to measure human progress.  

If the axioms and rules of inference are fixed, then this 
kind of progress
can be achieved mechanically by brute
force, which is not very interesting.  The interesting case is allowing 
\textbf{new} axioms and concepts.
So I would propose instead that
human progress---purely intellectual, not \textbf{moral} progress---be 
measured by the number of bits of $\Omega$ 
that we have been able to determine up to any given time.

Let me end with Leibniz's remarks about the effects of this kind of progress [33, pp.\ 115]:
\begin{quote}
If this happens, it must follow that those minds which are not yet sufficiently capable
will become more capable so that they can comprehend and invent such great theorems,
which are necessary to understand nature more deeply and to reduce physical truths to
mathematics, for example, to understand the mechanical functioning of animals, to forsee
certain future contingencies with a certain degree of accuracy, and to do certain
wonderful things in nature, which are now beyond our capacity\ldots

Every mind has a horizon in respect to its present intellectual capacity but not in
respect to its future intellectual capacity.
\end{quote}

\section*{References\footnote
{Many of the author's publications are available at his website at 
\begin{center}
\texttt{http://www.cs.auckland.ac.nz/CDMTCS/chaitin}
\end{center}
}}

\begin{itemize}
\item[{[1]}]
G. J. Chaitin, \emph{Exploring Randomness,} Springer-Verlag, 2001.
\item[{[2]}]
G. J. Chaitin, \emph{The Unknowable,} Springer-Verlag, 1999.
\item[{[3]}]
G. J. Chaitin, \emph{The Limits of Mathematics,} Springer-Verlag, 1998.
\item[{[4]}]
G. J. Chaitin, \emph{Conversations with a Mathematician,} Springer-Verlag, 2002.
\item[{[5]}]
G. J. Chaitin, ``On the intelligibility of the universe and the notions of 
simplicity, complexity and irreducibility,'' German Philosophy Association, in press.
\item[{[6]}]
M. du Sautoy, \emph{The Music of the Primes,} HarperCollins, 2003.
\item[{[7]}]
G. J. Chaitin, \emph{From Philosophy to Program Size,} Tallinn Cybernetics Institute, in press.
\item[{[8]}]
H. Weyl, \emph{The Open World,} Yale University Press, 1932, Ox Bow Press, 1989.
\item[{[9]}]
H. Weyl, \emph{Philosophy of Mathematics and Natural Science,}
Princeton University Press, 1949.
\item[{[10]}]
G. W. Leibniz, \emph{Philosophical Essays,} Hackett, 1989.
\item[{[11]}]
Leibniz, \emph{Discours de m\'etaphysique,} Gallimard, 1995.
\item[{[12]}]
M. Davis, \emph{The Universal Computer: The Road from Leibniz to Turing,} Norton, 2000.
\item[{[13]}]
J. L. Casti, \emph{The One True Platonic Heaven,} Joseph Henry Press, 2003.
\item[{[14]}]
E. T. Bell, \emph{Mathematics, Queen and Servant of Science,} Tempus, 1951.
\item[{[15]}]
G. J. Chaitin, \emph{Algorithmic Information Theory,} Cambridge University Press, 1987.
\item[{[16]}]
T. Ord and T. D. Kieu,
``On the existence of a new family of diophantine equations for $\Omega$,''
{\small \texttt{http://arxiv.org/abs/math.NT/0301274}}.
\item[{[17]}]
P. Raatikainen, Book review, \emph{AMS Notices} \textbf{48}, Oct.\ 2001, pp.\ 992--996.
\item[{[18]}]
H.-U. Obrist, \emph{HuO Interviews 1,} Charta, in press.
\item[{[19]}]
G. J. Chaitin, ``Meta-mathematics and the foundations of mathematics,''
\emph{EATCS Bulletin} \textbf{77}, June 2002, pp.\ 167--179.
\item[{[20]}]
O. Postel-Vinay, ``L'Univers est-il un calculateur?''\ [Is the universe a calculator?], 
\emph{La Recherche} \textbf{360}, Jan.\ 2003, pp.\ 33--44.
\item[{[21]}]
M. A. Nielsen and I. L. Chuang, \emph{Quantum Computation and Quantum Information,}
Cambridge University Press, 2000.
\item[{[22]}]
S. Wolfram, \emph{A New Kind of Science,} Wolfram Media, 2002.
\item[{[23]}]
L. Smolin, \emph{Three Roads to Quantum Gravity,} Weidenfeld and Nicolson, 2000.
\item[{[24]}]
T. Tymoczko, \emph{New Directions in the Philosophy of Mathematics,} Princeton University Press, 1998.
\item[{[25]}]
J. Borwein and D. Bailey, \emph{Mathematics by Experiment,} A. K. Peters, in press.
\item[{[26]}]
J. Borwein and D. Bailey, \emph{Experimentation in Mathematics}, A. K. Peters, in press.
\item[{[27]}]
D. S. Robertson, \emph{The New Renaissance,} Oxford University Press, 1998.
\item[{[28]}]
D. S. Robertson, \emph{Phase Change,} Oxford University Press, 2003. 
\item[{[29]}]
J. Maynard Smith and E. Szathm\'ary, \emph{The Major Transitions in Evolution,}
Oxford University Press, 1995.
\item[{[30]}]
J. Maynard Smith and E. Szathm\'ary, \emph{The Origins of Life,} Oxford University Press, 1999.
\item[{[31]}]
C. S. Calude, \emph{Information and Randomness,} Springer-Verlag, 2002.
\item[{[32]}]
E. T. Bell, \emph{Men of Mathematics,} Simon and Schuster, 1937.
\item[{[33]}]
A. P. Coudert, \emph{Leibniz and the Kabbalah,} Kluwer Academic, 1995.
\end{itemize}

\end{document}